\documentclass[12pt]{article}
\usepackage[dvips]{graphics}
\newcommand{\nc}{\newcommand}
\nc{\dl}{\delta} \nc{\ph}{\Phi_n} \nc{\la}{\lambda} \nc{\vp}{\varphi}
\nc{\bq}{\begin{equation}} \nc{\eq}{\end{equation}}
\nc{\ov}{\over} \nc{\hf}{{1\ov2}} \nc{\inv}{^{-1}} \nc{\tl}{\widetilde}
\renewcommand{\sp}{\vspace{1ex}} \nc{\iy}{\infty} \nc{\cd}{\cdots}
\nc{\La}{\Lambda} \nc{\noi}{\noindent}
\nc{\ba}{\begin{array}} \nc{\ea}{\end{array}} \nc{\pl}{\partial}
\nc{\T}{T_n(f)} \nc{\Ti}{T_n(f)\inv} \nc{\Tt}{T_n(\tl f)}
\nc{\Tti}{T_n(\tl f)\inv}
\nc{\fp}{f^+} \nc{\fm}{f^-} \nc{\ftp}{\tl f^+} \nc{\ftm}{\tl f^-}
\nc{\dlp}{\dl^+} \nc{\dlm}{\dl^-} \nc{\tn}{\otimes}
\nc{\Vp}{V^+} \nc{\Vm}{V^-} \nc{\Vtp}{\tl V^+} \nc{\Vtm}{\tl V^-}
\nc{\Up}{U^+} \nc{\Um}{U^-} \nc{\Utp}{\tl U^+} \nc{\Utm}{\tl U^-}
\nc{\twotwo}[4]{\left(\begin{array}{cc}#1&#2\\&\\#3&#4\end{array}\right)}
\nc{\twoone}[2]{\left(\begin{array}{c}#1\\\\#2\end{array}\right)}
\nc{\ra}{\rightarrow} \nc{\cE}{{\cal{E}}} \nc{\s}{\sigma}
\nc{\Dl}{\Delta}
\nc{\Ga}{\Gamma} \nc{\Z}{{\cal Z}}
\nc{\cS}{{\cal{S}}} \nc{\ga}{\gamma}
\nc{\cP}{{\cal{P}}} \nc{\ve}{\varepsilon} \nc{\om}{\omega}
\newcommand{\ban}{\begin{eqnarray*}}
\newcommand{\ean}{\end{eqnarray*}}
\begin{document}

\begin{center} {\large \bf On the Distributions of the Lengths of the
Longest\\\vspace{1ex}
Monotone Subsequences in Random Words }\end{center}\sp

\sp\begin{center}{{\bf Craig A.~Tracy}\\
{\it Department of Mathematics and Institute of Theoretical Dynamics\\
University of California, Davis, CA 95616, USA\\
e-mail address: tracy@itd.ucdavis.edu}}\end{center}
\begin{center}{{\bf Harold Widom}\\
{\it Department of Mathematics\\
University of California, Santa Cruz, CA 95064, USA\\
e-mail address: widom@math.ucsc.edu}}\end{center}\sp

\begin{center}{\bf Abstract}\end{center}
We consider the distributions of the lengths of the longest weakly
increasing and strongly decreasing subsequences in words of length
$N$ from an alphabet of $k$ letters. (In the limit as $k\ra\iy$
these become the corresponding distributions for permutations on $N$
letters.) We find Toeplitz determinant representations for
the exponential generating functions (on $N$) of these
distribution functions and show that they are expressible in terms
of solutions of Painlev\'e V equations. We show further that
in the weakly increasing case
the generating function gives the distribution of the
smallest eigenvalue in the $k\times k$ Laguerre random matrix
ensemble and that the
distribution itself has, after centering and normalizing, an
$N\ra\iy$ limit which is equal to the distribution function
for the largest eigenvalue in the
Gaussian Unitary Ensemble of  $k\times k$ hermitian matrices of trace
zero.
\sp

\begin{center} {\bf I. Introduction}\end{center}
\setcounter{equation}{0}\renewcommand{\theequation}{1.\arabic{equation}}

The last decade has seen a flurry of activity centering around
connections between combinatorial probability on the one hand and
random matrices and integrable systems on the other. From
the point of view of probability theory, the quite surprising
feature  of these developments is that the methods  
 came from  Toeplitz determinants,
integrable differential equations of the Painlev\'e type and the closely
related Riemann-Hilbert techniques as they were applied and
refined in random matrix
theory.  Using these methods new, and apparently quite universal, 
limiting laws have been discovered.  
One of the aims of this paper is
to make these methods accessible to a wider audience.  Our story
begins with a theorem of Gessel~\cite{gessel}. (There are
 earlier signs of these connections; see Regev~\cite{regev}.)

Let $S_N$ be the symmetric group on $N$ letters and give each
permutation
$\s\in S_N$ probability $1/N!$. Denote by $\ell_N(\s)$ the length of the
longest increasing subsequence in $\s$ and
\[ F_P(n;N)=\mbox{\rm Prob}\left(\ell_N(\sigma)\le n\right).\]
Then it is a corollary of Gessel's theorem and the
Robinson-Schensted-Knuth
correspondence\footnote{This is a bijection between permutations and
pairs $(P,Q)$ of  standard Young tableaux with the same shape. For
expository
accounts
of the RSK algorithm, see \cite{fulton,knuthBook,stanley}.} that
\bq \sum_{N=0}^\iy F_P(n;N) {t^N\ov N!}=D_n(t),\label{Pgenrep}\eq
where $D_n(t)$ is the determinant of the the $n\times n$ Toeplitz matrix
with the symbol
$e^{\sqrt t\,(z+z\inv)}$.
(Recall that the $i,\,j$ entry of a Toeplitz matrix equals the $i-j$
Fourier
coefficient of its symbol.)

It is in this work of Gessel, expressing the (exponential) generating
function
of $F_P$  as a Toeplitz determinant, and the subsequent work of
Odlyzko \textit{et al.}~\cite{odlyzko} and Rains~\cite{rains},
that the methods of random matrix theory first appear in RSK
type problems.\footnote{Gessel \cite{gessel}
does not mention random matrices,
but in light of well-known formulas in random matrix
theory relating Toeplitz determinants to expectations over the unitary
group, we
believe it is fair to say that the
connection with random matrix theory begins
with his discovery.}

Starting with this representation, Baik, Deift and Johansson
\cite{bdj1}, using the
steepest descent method for Riemann-Hilbert problems \cite{dz}, derived
a delicate asymptotic formula for $D_n(t)$ which we now describe.
Introduce another parameter
$s$ and suppose that $n$ and $t$ are related by $n=[2 t^{1/2} + s
t^{1/6}]$.
Then
as $t\ra\iy$ with $s$ fixed one has
\[ \lim_{t\ra\iy} e^{-t} D_{2\,t^{1/2}+s\,t^{1/6}}(t) = F_2(s). \]
Here $F_2$ is the distribution function defined by
\bq F_2(s) =\exp\left(-\int_s^\iy (x-s) q(x)^2\, dx\right)
\label{painF2}\eq
where $q$ is the solution of the Painlev\'e II equation
\[ q^{\prime\prime}=sq + 2 q^3 \]
satisfying $q(s)\sim \mbox{\rm Ai}(s)$ as $s\ra\iy$.\footnote{$\mbox{\rm
Ai}$
 is the Airy function.  For a proof that such a solution
exists and is unique, see~\cite{clarkson,hastings,dz2}.}  Using a
dePoissonization
lemma due to Johansson~\cite{johansson1}, these asymptotics for
$D_n(t)$ led Baik, Deift and Johansson  to the limiting law
\[ \lim_{N\ra\iy}\mbox{\rm Prob}\left({\ell_N(\sigma)-2\sqrt{N}\ov
N^{1/6}}\le s\right)
=F_2(s). \]

It is a remarkable fact that this same distribution function $F_2$ was
first encountered
by the present authors~\cite{tw1} in random matrix theory where it
arises as
the limiting law  for the normalized largest
eigenvalue in the Gaussian Unitary Ensemble (GUE) of Hermitian
matrices.  More precisely,
we have for this ensemble~\cite{tw1}
\bq \lim_{N\ra\iy}\mbox{\rm Prob}\left((\la_{\mbox{\scriptsize
max}}-\sqrt{2N})
\sqrt{2} N^{1/6}\le s\right)=F_2(s).\label{F2}\eq

Here we see a connection with integrable systems---the appearance of a
Painlev\'e II
function. Yet another connection is that $D_n(t)$ itself has a
representation in terms
of the solution of a Painlev\'e V equation~\cite{hisakado,tw4}.

Since the work of Baik, Deift and Johansson, several groups have
extended this connection
between RSK type combinatorics and the distribution functions of random
matrix theory.
The aforementioned result is equivalent to the determination of
the limiting distribution of the number of
boxes in the first row in the RSK correspondence $\sigma\leftrightarrow
(P,Q)$.
In \cite{bdj2} the same authors show that the limiting distribution of
the number of
boxes in the {\it second} row  is (when centered and normalized)
distributed as the
{\it second} largest scaled eigenvalue in GUE \cite{tw1}.
They then conjectured that this correspondence extends to all rows.
This
conjecture was recently proved by Okounkov \cite{okounkov} using
topological methods
and by Borodin, Okounkov and Olshanski \cite{borodin2} using
analytical/representation
theoretic methods.

Placing restrictions on the permutations $\sigma$ (that they be fixed
point free and
involutions), Baik and Rains~\cite{baikRains} have shown that the
limiting laws for
the length of the longest increasing/decreasing
subsequence are now the limiting distributions $F_1$ and $F_4$
\cite{tw3} for the
scaled largest eigenvalues in the Gaussian Orthogonal Ensemble (GOE) and
the Gaussian
Symplectic Ensemble (GSE).
Generalizing to signed permutations and colored permutations the present
authors and
Borodin \cite{tw4,borodine1} showed that the
distribution functions of the length of the longest
increasing
subsequence involve the same $F_2$. Johansson~\cite{johansson2} showed
that the shape
fluctuations of a certain random growth model, again appropriately
scaled,
converges in distribution to $F_2$.  (This random growth model is
intimately related to
certain randomly growing Young diagrams.)

Finally, we mention the work of Aldous and Diaconis \cite{aldous} where
they describe a
certain one-person card game, called ``patience sorting,''  which is
connected to these ideas
linking Young tableaux with the length of the longest increasing
subsequence in either
a random permutation or a random word; and thence to the limiting
distributions of largest eigenvalues.

At last we come to the subject of the present paper, which is the
question of what
can be said when instead of a random permutation on $N$ letters we have
a random
word of length $N$ from an alphabet of $k$ letters. This may be
thought of
as a function from $\{1,\,2,\cdots,N\}$ to $\{1,\,2,\cdots,k\}$ and it
is clear
what is meant by a (strictly or weakly) increasing or decreasing
subsequence. Unlike
the case of permutations there is a difference between the two. Given
such a word $w$ we
define
$\ell^I_N(w)$ to be the length of the longest {\it weakly increasing}
subsequence
in $w$, and define $\ell^D_N(w)$ to be the length of the longest {\it
strictly
decreasing} subsequence in $w$.\footnote{We could just as
well consider weakly increasing subsequences \textit{and}
strictly increasing subsequences since reversing the order
of the word takes increasing subsequences into decreasing
subsequences.  In light of the RSK algorithm, our choice
seems most convenient.}
 In analogy with permutations we define
the
distribution functions (giving each word probability $k^{-N}$)
\[F_I(n;k,N)=\mbox{\rm Prob}_k\left(\ell_N^I(w)\le n\right),\ \ \ \
F_D(n;k,N)=\mbox{\rm Prob}_k\left(\ell_N^D(w)\le n\right)\]
and their generating functions
\bq G_I(n;k,t)=\sum_{N=0}^{\iy}F_I(n;k,N)\,{t^N\ov N!},\ \ \ \ \
G_D(n;k,t)=\sum_{N=0}^{\iy}F_D(n;k,N)\,{t^N\ov N!}.\label{G}\eq
Here are our results. We use the standard notation $T_n(f)$ for the
$n\times n$ Toeplitz matrix with symbol~$f$.\sp

\noi{\bf Theorem 1}. We have
\[G_I(n;k,kt)=\det\,T_n(f_I),\ \ \ G_D(n;k,kt)=\det\,T_n(f_D),\]
where
\[f_I(z)=f_I(z;t)=e^{t/z}\,(1+z)^k,\ \ \
f_D(z)=f_D(z;t)=e^{t/z}\,(1-z)^{-k}.\]\sp

\noi{\bf Theorem 2}. $G_I(n;k,t)$ and $G_D(n;k,t)$ have integral
representations in terms of solutions of Painlev\'e V equations.\sp

\noi{\bf Theorem 3}. $G_I(n;k,t)$ is equal to $e^{kt}$ times the
distribution function
for the smallest eigenvalue in the Laguerre ensemble of $k\times k$
matrices associated with the weight function $x^n\,e^{-x}$.\sp

\noi{\bf Theorem 4}. The limiting distribution for the random variable
$\ell^I_N(w)$,
centered and normalized, is equal to that for the largest eigenvalue in
the Gaussian Unitary Ensemble ensemble of $k\times k$ hermitian matrices 
with trace zero.\footnote{The
normalization we adopt for the GUE measure is the standard
one defined in Mehta~\cite{mehta}. The probability on the right of the
displayed formula is the conditional probability given that the matrix
from GUE has trace
zero.} More precisely,
\[ \lim_{N\ra\iy}\mbox{Prob}_k\left({\ell_N^I(w)-N/k\ov \sqrt{2N/k}}\le
s\right)
=\mbox{Prob}\left(\la_{\mbox{\scriptsize max}}\le s\right).\]\sp

The next four sections contain the proofs of these four theorems.

Theorem 1 is a consequence of Gessel's theorem, just as the permutation
analogue (\ref{Pgenrep}) is, and the RSK correspondence between words
and pairs of Young tableaux. For the convenience
of the reader we include a complete proof of Gessel's theorem,
containing the main
ideas of the original but presented somewhat differently.
(For related developments, see \cite{baikRains2}.)

Theorem 2 is the heart of the paper. The equations are derived very much
in
the spirit of \cite{tw4}. The logarithmic derivative of the determinant
involves
a quantity whose derivatives in turn involve other quantities. Recursion
formulas relating the various quantities allow the eventual derivation
of a single
differential equation which, in the end, turns out to be reducible to
Painlev\'e V.

The proof of Theorem 3 consists of showing that the P$_V$ function
of Theorem 2 for $G_I$ is exactly the one which gives the distribution
of
the smallest eigenvalue in the Laguerre ensemble \cite{tw2}. The
equation is the
same, by inspection, and it is a matter of checking the boundary
condition at $t=0$.

Given the results of~\cite{bdj2, okounkov, borodin2} for
permutations, it is natural to conjecture that
the limiting distribution of the
number of boxes in the $j^{th}$ row, $2\le j \le k$, appearing
in the Young tableaux $P$ in the RSK bijection $w\leftrightarrow (P,Q)$
is precisely the distribution of the $j^{th}$ largest eigenvalue in
the finite $k\times k$ Hermite ensemble.

Theorem 4 is proved by an asymptotic evaluation of the multiple
integral which gives the distribution function for the smallest
eigenvalue in the Laguerre ensemble.

After the completion of the original version of this paper there were
several relevant
developments. Johansson \cite{johansson3} found an independent proof of
Theorem 4, in fact of the full conjecture stated above. A.~Its, using
Riemann-Hilbert
techniques applied to operator equivalents of our Toeplitz matrices,
found another
proof of Theorem 2. J.~L.~Snell found an error in the
original version of Theorem 4. (We thank him for catching the error and
so saving the authors
from further embarrassment.) C.~Grinstead found a random walk
interpretation
of the $k=2$ problem and used this to determine the limiting
distribution
in this case. And Y.~Chen alerted us to the paper \cite{forrester} of
Forrester in which
there appeared a formula equivalent (given Theorem~1) to the statement
of Theorem 3, 
obtained by using identities of Macdonald on hypergeometric
functions of several variables.

\sp

\begin{center} {\bf II. Gessel's theorem and its
specializations}\end{center}
\setcounter{equation}{0}\renewcommand{\theequation}{2.\arabic{equation}}

\begin{center}{\bf 1. The Cauchy-Binet formula}\end{center}

We begin by recalling
the Cauchy-Binet formula for the determinant of
the product of two rectangular matrices $A$ and
$B$ of sizes $m\times n$ and $n\times m$, respectively.  We assume $n\ge
m$.

Let $\cS_{mn}$ denote the set of strictly increasing subsequences
of length $m$ that can be chosen from $\{1,2,\ldots,n\}$.  For any
matrix $X$
of size $n\times m$ and any $S=\{s_1,s_2,\ldots,s_m\}\in\cS_{mn}$,
denote by $X(S\vert m)$ the $m\times m$ matrix obtained from $X$ by
using all $m$ columns
of $X$ and the $m$ rows numbered by $S$.  Similarly, if $X$ is $m\times
n$, denote
by $X(m\vert S)$ the $m\times m$ matrix obtained from all $m$ rows of
$X$ and the
columns of $X$ numbered by $S$.  The Cauchy-Binet formula is
\bq \det(AB)=\sum_{S\in\cS_{mn}}\det\left(B(S\vert m)\right)
 \det\left(A(m\vert S) \right).
\label{cb}\eq

Here is a proof. Define
\[ F(\ve)=\ve^m\det\left(I+\ve^{-1}AB\right).\]
This determinant is $m\times m$ and $\lim_{\ve\ra 0}F(\ve)=\det(AB)$.
Then also
\[ F(\ve)=\ve^m\det\left(I+\ve^{-1} BA\right). \]
In the (Fredholm) expansion of this last determinant, the term with
coefficient $\ve^{-m}$
is
\[ {1\ov m!}\sum_{i_1,i_2,\ldots,i_m}\left| \begin{array}{llcl}
        (BA)_{i_1 i_1} & (BA)_{i_1 i_2} & \cdots & (BA)_{i_1 i_m} \\
        (BA)_{i_2 i_1} & (BA)_{i_2 i_2} & \cdots & (BA)_{i_2 i_m} \\
        \multicolumn{4}{c} \dotfill \\
        (BA)_{i_m i_1} & (BA)_{i_m i_2} & \cdots & (BA){i_m i_m}
        \end{array}\right| \]
In this sum we can place the restriction that no two indices are equal
since when they are the determinant is zero.  The $m!$ different
orderings
of $\{i_1,\ldots,i_m\}$ give the same determinant so we can drop
the $m!$ and sum over all $i_\alpha$ with $i_1<i_2<\cdots<i_m$.
The terms of order $\ve^{-j}$, $j<m$, in the Fredholm expansion do not
contribute
to the limit as $\ve\ra 0$.  The coefficients of the terms $\ve^{-j}$,
$j>m$, are
zero since the rank of $BA$ is at most $m$.  Finally, each summand
factors into the product of the two determinants in the Cauchy-Binet
formula.

The formula remains valid if $n=\infty$ if, for example, each
row of $A$ and each column of $B$ belongs to the sequence space
$\ell^2$. For then
$AB$ is well-defined, $BA$ is a finite rank operator on $\ell^2$ and the
preceding goes
through without change. The sum on the
right side of (\ref{cb}) is then the sum over all increasing
subsequences $S$ of
length $m$ chosen from the positive integers.

\begin{center}{\bf 2. Gessel's theorem}\end{center}

Let $\cP_m$ denote the set of  partitions of $m$, sequences
$(\la_1,\,\la_2,\,\ldots,\la_n)$
of nonnegative integers such that $\la_1\ge\cdots\ge\la_n,\
\la_1+\cdots+\la_n=m$,
and $\cP=\bigcup_{m\ge 0} \cP_m$.  ($\cP_0$ is the empty partition.)
We also write $\la\vdash m$ when $\la\in\cP_m$ and denote by
$\ell(\lambda)$ the length of the partition, the largest $k$ such that
$\la_k\ne0$.
Let $\Lambda_{\mbox{\bf\footnotesize{Q}}}$ (or $\Lambda$ for short)
denote the
algebra of symmetric functions over {\bf Q}.  This is a commutative
algebra and the vector space direct sum decomposition into
homogeneous symmetric functions gives
$\Lambda_{\mbox{\bf\footnotesize{Q}}}$ the
structure of a graded algebra.

Gessel introduces
\bq R_n(x,y):=\sum_{{\lambda\in\cP \atop \ell(\lambda)\le n}}
s_\lambda(x) s_\lambda(y)=\sum_{m=0}^\iy
\sum_{{\lambda\in\cP_m \atop \ell(\lambda)\le n}}
s_\lambda(x) s_\lambda(y)
 \label{R}\eq
where $s_\lambda(x)$ is the Schur function.  Gessel's theorem says
that $R_n(x,y)$ is a Toeplitz determinant
\bq R_n(x,y)=\det\left(A_{i-j}\right)_{1\le i,j\le n}\label{gessel}\eq
where
\[ A_i=A_i(x,y)=\sum_{\ell=0}^\iy h_{\ell+i}(x) h_{\ell}(y) \]
and $h_r$ is the $r^{\mbox{\rm\footnotesize{th}}}$ complete symmetric
function.
(We take $h_r=0$ for $r<0$.)  Recall
that
\[ \sum_{r\ge 0} h_r t^r = \prod_{i\ge 1} \left(1-x_i t\right)^{-1}.\]

Here is Gessel's proof. Although (\ref{gessel}) is an identity between
two formal power
series it suffices to prove it when the $x_i$ are real numbers
satisfying $|x_i|<1$.
Let $M(x)$ be the $\iy\times n$ matrix
\[(h_{i-j}(x)),\ \ \ \ (i\ge 1,\ 1\le j \le n).\]
This will be the matrix $B$ of (\ref{cb}) whereas $A$ will be $M^t(y)$.
(We interchanged
the roles of $m$ and $n$.) Since the entries of the columns of $M(x)$
and
rows of $M^t(y)$
are exponentially small, the formula holds.

For any increasing subsequence $S$ of positive
integers of length $n$, let $M_S(x)$ be the determinant of the $n\times
n$ minor
of $M(x)$ obtained from the rows indexed by the elements of $S$. In the
notation of (\ref{cb}) we have $M_S(x)=\det(B(S\vert n)),\
M_S(y)=\det(A(n\vert S))$.

Now let $\lambda=(\lambda_1,\ldots,\lambda_n)$ be a partition with
$\ell(\lambda)\le n$, and let
$S=\{\lambda_{n+1-i}+i\vert 1\le i \le n\}$. Observe that $S$
is an increasing subsequence of postive integers.  Then
\[M_S(x)=\det\left(h_{\lambda_{n+1-i}+i-j}(x)\right)_{1\le i,j\le n}.\]
Reversing
the order of the rows and columns in this determinant yields
\[M_S(x)=\det\left(h_{\lambda_i+j-i}(x)\right)=s_\lambda(x),\]
where the last equality is the Jacobi-Trudi identity.
For such $S$,
summing over all partitions $\lambda$ with $\ell(\lambda)\le n$ is the
same as summing over all increasing subsequences $S$ of length $n$.
Thus
\[ R_n(x,y)=\sum_{S} M_S(x) M_S(y), \]
and by the Cauchy-Binet formula this is equal to $\det(M^t(y)\, M(x))$.
Since
\[ (M^t(y)\, M(x))_{ij}=\sum_{\ell} h_{\ell-i}(y)\,h_{\ell-j}(x)
=\sum_{\ell} h_{\ell+i-j}(x)\,h_{\ell}(y)=A_{i-j},\]
the theorem follows.

For the applications which follow it is important to know the the symbol
of
the Toeplitz determinant appearing in (\ref{gessel}). It is
\[\varphi(z)=\sum_{i=-\iy}^\iy A_i(x,y)\, z^i
=\sum_{i=-\iy}^\iy z^i \sum_{\ell=0}^\iy h_{\ell+i}(x) h_\ell(y)
=\sum_{\ell=0}^\iy h_{\ell}(y)\, \sum_{i=-\iy}^\iy h_i(x) z^{i-\ell}\]
\bq=\prod_{n=1}^\iy \,\left(1-y_n\,z^{-1}\right)^{-1}
\prod_{n=1}^\iy \left(1-x_n\,z\right)^{-1}.\label{symbol}\eq\sp

\begin{center}{\bf 3. Cauchy's identity from Szeg\"o's
theorem}\end{center}

A nice application of Gessel's theorem is a derivation of Cauchy's
identity
in symmetric functions\footnote{We
freely use various results from symmetric functions which
can be found, for example, in~\cite{stanley}, Chp.~7.}
 using the strong Szeg\"o's limit theorem for
Toeplitz determinants.
We have $(\log\varphi)_0=0$ and, for $n>0$,
\[(\log\varphi)_n={1\over n} \sum_{i\ge 1} (x_i)^n,\]
\[(\log\varphi)_{-n}={1\ov n}\sum_{i\ge 1} (y_i)^n.\]
(The subscripts denote Fourier coefficients, as usual.)
Applying Szeg\"o's theorem (we may assume that the $x_i$ and $y_i$ are
real numbers with absolute value less than 1) then gives
\[ \lim_{n\ra\iy} R_n(x,y)=\exp\left\{\sum_{n=1}^{\infty}{1\ov
n}\sum_{i,\,j\ge1}
(x_i\,y_j)^n\right\}=\prod_{i,\,j}\left(1-x_i y_j\right)^{-1},\]
and hence  Cauchy's identity
\[\sum_{\lambda\in\cP}
s_\lambda(x) s_\lambda(y) =\prod_{i,\,j}
\left(1-x_i y_j\right)^{-1}.\]\sp\pagebreak

\begin{center}{\bf 4. Dual version of Gessel's theorem}\end{center}

Since (\ref{gessel}) is an identity between two elements of the ring
generated by
the $x_i$ and $y_i$, any
endomorphism of this ring yields another identity. Now the complete
symmetric
functions $h_r$ are algebraically independent generators of $\La$ as are
the elementary
symmetric functions $e_r$. We consider the endomorphism $\om$ defined by
\[\om(e_r)=h_r.\]
Then for any partition $\la=(\la_1,\,\la_2,\ldots)$ we have
\[\om(e_{\la})=h_{\la}\]
with the usual notation
\[e_{\la}=\prod_ie_{\la_i},\ \ h_{\la}=\prod_ih_{\la_i}.\]
The action on the Schur function is given by $\om(s_\la)=s_{\la'}$ where
$\la'$ is
the partition conjugate to $\la$.

We define
\[\tilde{R}_n(x,y)=\om_x \, \om_y \, R_n(x,y)
=\sum_{{\la\in\cP \atop \ell(\lambda) \le n}}
s_{\la'}(x) s_{\la'}(y).\]
Now for any partition $\la$ we have $ \ell(\lambda)=\la'_1$, the length
of the first
row of the Young diagram of shape $\la'$, and so
\[\tilde{R}_n(x,y)=\sum_{{\la\in\cP \atop \la'_1\le n}}s_{\la'}(x)
s_{\la'}(y)
=\sum_{{\la\in\cP \atop \la_1\le n}}s_{\la}(x) s_{\la}(y).\]
(Thus, in this sum we restrict the length of the first row rather than
the
length of the first column.) Applying $\om_x\,\om_y$ to (\ref{gessel})
we obtain
the dual version of Gessel's theorem,
\[\sum_{{\la\in\cP \atop \la_1\le n}}s_{\la}(x) s_{\la}(y)=
\det\,\left(\om_x\,\om_y\,A_{i-j}(x,y)\right)_{i,\,j=1,\ldots,n}.\]
We record for use below
\[\om_x \, \om_y \,\vp(z)=
\sum_{\ell=0}^\iy e_{\ell}(y)\, \sum_{i=-\iy}^\iy e_i(x) z^{i-\ell}\]
\bq=\prod_{n=1}^\iy \,\left(1+y_n\,z^{-1}\right)
\prod_{n=1}^\iy \left(1+x_n\,z\right).\label{dualsymbol}\eq

We remark that $\omega_y R_n(x,y)$
also equals a Toeplitz determinant.  The $n\ra\iy$ limit of
this identity is,  by
an application of the strong Szeg\"o theorem,  the \textit{dual}
Cauchy identity
\[ \sum_{\lambda\in\cP} s_\lambda(x) s_{\lambda^\prime}(y)=
\prod_{i,j}(1+x_i y_j). \]
\sp

\begin{center}{\bf 5. Specializations}\end{center}

If $R$ is a commutative {\bf Q}-algebra
with identity, then  a {\it specialization} of the ring $\Lambda$
is a homomorphism $\psi:\Lambda\ra R$.  We always assume
that $\psi(1)=1$.\sp

\noi{\it (i) Exponential specialization}.  If $p_n=\sum_i x_i^n$, the
power sum symmetric functions, then the exponential specialization is
determined by
\[ ex(p_n)= t \delta_{1n}. \]
(Recall that the $p_n$'s are algebraically independent generators of
$\La$.
This homomorphism $ex$ is denoted by $\theta$ in Gessel.)

The fundamental
property of this homomorphism is for any symmetric function $f$
\[ ex(f) = \sum_{n\ge 0} [x_1x_2\cdots x_n]f\, {t^n\ov n!} \]
where $[x_1x_2\cdots x_n]f$ denotes the coefficient of $x_1x_2\cdots
x_n$
in $f$.  Thus if $\lambda\vdash N$ and $s_\lambda$ is the Schur
function, then
\[ ex(s_\lambda)= f^\lambda\, {t^N\ov N!} \]
where $f^\lambda$ is the number of standard Young tableaux of
shape $\lambda$. Hence
\[ ex_x\, ex_y R_n(x,y)=\sum_{N=0}^\iy \sum_{{\la\vdash N
\atop \ell(\lambda)\le n}} \left(f^\lambda\right)^2 {t^{2N}\ov
(N!)^2}.\]

By the Robinson-Schensted-Knuth (RSK)
bijection~\cite{robinson,schensted,knuth}
\[\sum_{{\la\vdash N\atop \ell(\lambda)\le n}} \left(f^\la\right)^2\]
equals the number of permutations $\s$ on $N$ letters such that
$\ell_N(\s)$, the
length of the longest increasing subsequence in $\s$, is at most $n$.
Hence if
each such permutation has probability $1/N!$ we have
\[F_P(n;N):=\mbox{Prob}\left(\ell_N(\sigma)\le n\right)={1\ov N!}
\sum_{{\la\vdash N
\atop \ell(\la)\le n}} \left(f^\la\right)^2.\]
Thus we know that its generating function is given by
\[\sum_{N=0}^\iy F_P(n;N)\,{t^{2N}\ov N!}=ex_x\, ex_y R_n(x,y).\]

Gessel's theorem tells us that $R_n(x,y)$ is the $n\times n$ Toeplitz
determinant
with symbol $\vp(z)$ given by (\ref{symbol}). This may be written
\[\vp(z)=\sum_{r\ge 0} h_r(y) z^{-r}\, \sum_{s\ge 0} h_s(x) z^s.\]

The important observation is that since $ex$ is a homomorphism, $ex_x\,
ex_y R_n(x,y)$
is the Toeplitz determinant with symbol
\[f_P(z):=
ex_x\, ex_y (\varphi(z))=\sum_{r\ge 0} {t^r\ov r!} z^{-r}\,\sum_{s\ge 0}
{t^s\ov s!} z^s
= e^{t\,(z+z\inv)}.\]
This is precisely (\ref{Pgenrep}) after changing $t$ to $\sqrt t$.\sp

\noi{\it (ii) Principal specializations}.
The {\it principal specialization} of order $n$ of $g$ is defined by
\[ ps_n(g)=g(1,q,q^2,\ldots,q^{n-1},0,0,\ldots). \]
(Thus we replace $x_i$ by $q^i$ if $i<n$ and by 0 itherwise.
If we let $n\ra\iy$ we obtain the {\it stable principal
specialization\/}
of $f$.)
Setting $q=1$ in $ps_n$ gives
\[ ps_n^1(g)=g(1,1,\ldots,1,0,0,\ldots)\]
where 1 appears $n$ times.

Observe that
\[ ps_k^1(s_\lambda) = d_\lambda(k)\]
where $d_\lambda(k)$ is the number of semistandard Young
tableaux of shape $\lambda$ that can be formed  from an alphabet
of $k$ letters.  (Recall that a semistandard tableau is weakly
increasing
across rows and strictly increasing down columns; a standard tableau is
strictly
increasing across rows.) This is most easily seen from the combinatorial
definition of the Schur functions.

Applying the homomorphisms $ps_n^1$ and $ex$ to $R_n$ gives
\bq (ps_k^1)_x \, ex_y \,R_n(x,y)=
\sum_{N=0}^\iy \left(\sum_{{\lambda\vdash N
\atop \ell(\lambda)\le n}} d_\lambda(k)\, f^\lambda \right){t^N\ov
N!}.\label{gen1}\eq
The RSK correspondence associates to each word $w$ of length $N$
formed from an alphabet of $k$ letters a pair of tableaux, $(P,Q)$.
Here the $P$ are semistandard Young tableaux of shape $\lambda\vdash N$
made from the alphabet $\{1,2,\ldots,k\}$  and the $Q$ are
standard Young tableaux of shape $\lambda\vdash N$ on
the numbers $\{1,2,\ldots, N\}$. Thus
\[ \sum_{{\lambda\vdash N
\atop \ell(\lambda)\le n}} d_\lambda(k)\, f^\lambda \]
counts the number of words $w$ of length $N$ with strictly decreasing
subsequences all of length less than or equal to $n$.  Obviously,
\[ \sum_{\lambda\vdash N} d_\lambda(k)\, f^\lambda = k^N.\]

The symbol of the Toeplitz determinant that equals the generating
function (\ref{gen1}) is
\[f_D(z)=f_D(z;\,t):=(ps_k^1)_x ex_y(\varphi(z))=\prod_{i=1}^k
(1-z)^{-1}\,
\sum_{r\ge 0} {t^r\ov r!} z^{-r}={e^{t/z}\ov (1-z)^k}.\]
Hence we have shown that if $\ell_N^D(w)$ denotes the length
of the longest {\it strictly decreasing} subsequence in word $w$, and if
each
word of length $N$ is assigned probability $1/k^N$, then the
generating function of the distribution function
\[ F_D(n;k,N):=\mbox{Prob}_k\left(\ell_N^D(w)\le n\right)\]
is given by the Toeplitz determinant with symbol $f_D$:
\[\sum_{N=0}^\iy F_D(n;k,N)\,{(kt)^N\ov N!}
=\det\left(T_n(f_D(z;\,t)\right).\]
Recalling (\ref{G}) we see that this is
\bq G_D(n;k,kt)=\det\left(T_n(f_D(z;\,t)\right).\label{toeplitzD}\eq

Since, under general conditions, changing the symbol of a Toeplitz
matrix from $f(z)$
to $f(az)$ is a similarity
transformation, the associated Toeplitz determinant does not change.
Therefore the symbol of the Toeplitz determinant in (\ref{toeplitzD})
may be replaced by
 \[f(\sqrt{t}z/k;\,t/k)={e^{\sqrt{t}/z}\ov (1-\sqrt{t}z/k)^k } \]
 whose $k\ra\iy$ limit is $e^{\sqrt{t}\,(z+z\inv)}$.  This shows that
for fixed $N$,
 \[ \lim_{k\ra\iy} F_D(n;k,N)=F_P(n;N).\]
 Again, this is intuitively clear since as the size of the alphabet
approaches
 infinity, any random word of length $N$ is very likely a permutation.
 (This also uses the fact that the distribution of the length
 of the longest decreasing subsequence of a random permutation is the
same
 as the distribution of the length of the longest increasing
subsequence.)

Finally we apply the same specialization $(ps_k^1)_x \, ex_y$ to the
dual version
of Gessel's theorem. We see that (\ref{gen1}) is replaced by
\bq (ps_k^1)_x\, ex_y\, \tilde{R}_n(x,y)
=\sum_{N=0}^\iy \left(\sum_{{\lambda_1\le n \atop \lambda\in\cP}}
d_\lambda(k) f^\lambda\right) {t^N\ov N!}. \label{gen2}\eq
(We used here the fact that $ f^{\lambda^\prime}= f^\lambda$,
which follows immediately from the hook length formula for $f^\lambda$.)
Thus we obtain the generating function for the distribution
of the length of the longest {\it weakly increasing} subsequence.
Using (\ref{dualsymbol}) we find that the symbol of the Toeplitz
determinant that
gives (\ref{gen2}) is
\[f_I(z)=f_I(z;\,t)
:=(ps^1_k)_x\, ex_y \, \sum_{r\ge 0} e_r(x) z^r
\, \sum_{s\ge 0} e_s(y) z^{-s}\]
\[=(ps^1_k)_x\, ex_x \, \prod_{j}\left(1+x_j z\right)
\, \sum_{s\ge 0} e_s(y) z^{-s}
=(1+z)^k \, e^{t/z}.\]

To summarize, we have shown that if $\ell_N^I(w)$
denotes the length of the longest  weakly increasing
subsequence in word $w$ of length $N$, and if each
such word has probability $1/k^N$, then the generating
function of the distribution function
\[ F_I(n;k,N):=\mbox{Prob}_k\left(\ell_N^I(w)\le
n\right)\]
is given by the Toeplitz determinant with symbol $f_I$. Precisely,
\bq G_I(n;k,kt)=\sum_{N=0}^\iy F_I(n;k,N)\, {(kt)^N\ov N!}
=\det\left(T_n(f_I(z;\,t)\right).\label{toeplitzI}\eq
The same $k\ra\iy$ remarks hold here as in the strictly decreasing
case.

Relations (\ref{toeplitzD}) and (\ref{toeplitzI}) are the assertions of
Theorem 1.\sp

\begin{center} {\bf III. Recursion and differentiation
formulas}\end{center}
\setcounter{equation}{0}\renewcommand{\theequation}{3.\arabic{equation}}

\begin{center}{\bf 1. Universal recursion relations}\end{center}

In this section $f$ will be an arbitrary function, with Fourier
coefficients
$f_i$ and associated $n\times n$ Toeplitz matrix
\[\T=(f_{i-j}),\ \ \ \ \ (i,\,j=0,\cd,n-1).\]
We assume $\T$ is invertible and obtain several relations connecting
various
inner products involving $\Ti$. Most of these relations actually
appeared in~\cite{tw4}. There our $\T$ was symmetric and unfortunately
some of the relations derived in \cite{tw4} used this fact. Since this
does
not happen here,
everything has to be modified for the more general case. A reader of the
earlier
article will find familiar much of what we now do.

We introduce the $n$-vectors
\[\dlp=\left(\ba{c}1\\0\\\vdots\\0\\0\ea\right),\ \ \
\dlm=\left(\ba{c}0\\0\\\vdots\\0\\1\ea\right),
\ \ \ \fp=\left(\ba{c}f_1\\f_2\\\vdots\\f_{n-1}\\f_n\ea\right),\ \ \
\fm=\left(\ba{c}f_n\\f_{n-1}\\\vdots\\f_2\\f_1\ea\right)\]
and define $\tl f$ by $\tl f(z)=f(z\inv)$, so that $\Tt$ is the
transpose of $\T$.
We write
\[\La=T_n(z\inv),\ \ \ \La'=T_n(z).\]
Thus $\La$ is the backward shift and $\La'$ is the forward shift. It is
easy to see that
\bq T_n(z\inv\,f)=\T\,\La+\fp\tn\dl^+=\La\,T_n(f)+\dl^-\tn
f^-,\label{Lcom}\eq
\bq
T_n(z\,f)=T_n(f)\,\La'+\ftm\tn\dl^-=\La'\,T_n(f)+\dl^+\tn\ftp.\label{L'com}\eq
These identities explain why the vectors $f^{\pm}$ and
$\tl f^{\pm}$ arise.

The inner products involving $\Ti$ are
\[U_n^{\pm}=(\Ti \fp,\,\dl^{\pm}), \ \ \
V_n^{\pm}=(\Ti \dlp,\,\dl^{\pm}).\]
If one of these quantities defined in terms of $f$ is given a tilde,
then $f$
is to be replaced by $\tl f$ everywhere in its definition. Thus, for
example,
\[\Utp_n=(\Tti \ftp,\,\dlp).\]
Note that $\Vtp_n=\Vp_n=D_n/D_{n+1}$, where
\[D_n=D_n(f)=\det\,\T.\]
Some other inner products may be expressed in terms of these using the
isometry that
reverses the order of the components of a vector (and replaces a
Toeplitz matrix by
its transpose). Thus, for example,
\bq(\Tti \fm,\,\dl^+)=(\Ti \fp,\,\dlm)=\Um_n.\label{other}\eq
We shall use this isometry from time to time below without comment.

The basis for all the universal relations we shall obtain is the
following formula
for the inverse of a $2\times 2$ block matrix:
\bq\twotwo{A}{B}{C}{D}\inv=\twotwo{(A-BD\inv C)\inv}{\times}
{\times}{\times}.\label{ABCD}\eq
Here we assume $A$ and $D$ are square and the various inverses exist.
Only one block of the inverse is displayed and the formula shows that
$A-BD\inv C$
equals the inverse of this block of the inverse matrix.

We apply (\ref{ABCD}) first to the $(n+1)\times(n+1)$ matrix
\[\left(\begin{array}{cccc}0&0&\cd&1\\f_1&f_0&\cd&f_{-n+1}\\
\vdots&\vdots&\cd&\vdots\\f_n&f_{n-1}&\cd&f_0\end{array}\right),\]
with $A=(0),\ \ D=\T,\ \ B=(0\;\cd\;0\;1),\ \ C=f^+$.
In this case $A-BD\inv C=-(\Ti\,f^+,\,\dl^-)=-\Um_n$. This equals the
reciprocal of the upper-left entry
of the inverse matrix, which in turn equals $(-1)^n$ times the
lower-left $n\times
n$ subdeterminant divided by $D_n$. Replacing the first
row by $(f_0\;f_{-1}\;\cd\;f_{-n})$ gives the matrix
\bq\left(\begin{array}{cccc}f_0&f_{-1}&\cd&f_{-n}\\f_1&f_0&\cd&f_{-n+1}\\
\vdots&\vdots&\cd&\vdots\\f_n&f_{n-1}&\cd&f_0\end{array}
\right)=T_{n+1}(f).\label{Tn+1}\eq
The lower-left entry of its inverse equals on the one hand
$(\Ti\dlp,\,\dlm)=V_{n+1}^-$ and on the other
hand $(-1)^n$ times the same subdeterminant as arose above divided by
$\Um_{n+1}$. This gives the identity
\bq -\Um_n=V_{n+1}^-\,{D_{n+1}\ov D_n}={V_{n+1}^-\ov
V_{n+1}^+}.\label{UV}\eq

If we take $A$ to be the upper-left corner of (\ref{Tn+1}) and $D$
the complementary
$\T$ then $C=f^+$ and $B=(f_{-1}\;\cd\;f_{-n})$, and we deduce that
\bq f_0-(\Ti\fp,\,\ftp)={1\ov V_{n+1}^+}.\label{f0}\eq

Next we consider
\[\left(\begin{array}{ccccc}f_0&f_{-1}&\cd&f_{-n}&f_{-n-1}\\
f_1&f_0&\cd&f_{-n+1}&f_{-n}\\
\vdots&\vdots&\cd&\vdots&\vdots\\
f_n&f_{n-1}&\cd&f_0&f_{-1}\\
f_{n+1}&f_n&\cd&f_1&f_0\end{array}\right)=T_{n+2}(f).\]
We apply to this an obvious modification of (\ref{ABCD}), where $A$ is
the $2\times 2$ matrix consisting of the four corners of the large
matrix, $D$ is the
central $\T$, $C$ consists of the two columns $\fp$ and $\ftm$ and $B$
consists of the rows which are the transposes of $\ftp$ and $\fm$. Then
we find
\[A-BD\inv
C=\twotwo{f_0-(\Ti\fp,\,\ftp)}{f_{-n-1}-(\Ti\ftm,\,\ftp)}
{f_{n+1}-(\Ti\fp,\,\fm)}{f_0-(\Ti\ftm,\,\fm)}\]
and our formula tells us that this is the inverse of
\[\twotwo{V_{n+2}^+}{\tl V_{n+2}^-}{V_{n+2}^-}{V_{n+2}^+}.\]
This gives the two formulas
\[f_0-(\Ti\fp,\,\ftp)={V_{n+2}^+\ov {V_{n+2}^+}^2-V_{n+2}^-\tl
V_{n+2}^-}\,,\ \ \
f_{n+1}-(\Ti\fp,\,\fm)={-V_{n+2}^-\ov {V_{n+2}^+}^2-V_{n+2}^-\tl
V_{n+2}^-}\,.\]
Comparing the first with (\ref{f0}) we see that
\bq{V_{n+2}^+}^2-V_{n+2}^-\tl
V_{n+2}^-=V_{n+1}^+\,V_{n+2}^+,\label{V}\eq
and therefore that the preceding relations can be written
\bq f_0-(\Ti\fp,\,\ftp)={1\ov V_{n+1}^+},\ \ \
f_{n+1}-(\Ti\fp,\,\fm)=-{1\ov V_{n+1}^+}\,{V_{n+2}^-\ov
V_{n+2}^+}.\label{fn}\eq
Notice that (\ref{UV}) and (\ref{V}) give
\bq 1-\Um_n\,\Utm_n={V_n^+\ov V_{n+1}^+}.\label{UV1}\eq

Next, we apply (\ref{ABCD}) to the matrix
\[\left(\begin{array}{cccc}0&1&\cd&0\\f_1&f_0&\cd&f_{-n+1}\\
\vdots&\vdots&\cd&\vdots\\f_{n}&f_n&\cd&f_0\end{array}\right).\]
Now we have $A=(0),\ \ D=\T,\ B=(1\;\cd\;0\;0),\ \ C=f^+$
and so $A-BD\inv C=-(\Ti\,f^+,\,\dl^+)=-\Up_n$. Therefore this equals
the inverse of
the 0,\,0 entry of the inverse of
the matrix, which in turn equals its determinant divided by $D_n$. But
its determinant
equals $-K$, where $K$ is the 0,\,1 cofactor. Thus
\[U_n^+={K\ov D_n}.\]
Now look at the matrix (\ref{Tn+1})
and consider the 1,\,0 entry of its inverse. It equals on the one hand
$(T_{n+1}(f))\inv\,\dl^+,\,\La'\dl^+)$ and on the other hand
$-K/D_{n+1}$. This gives the
identity
\bq U_n^+=-{D_{n+1}\ov D_n}\,(T_{n+1}(f)\inv\,\dl^+,\,\La'\dl^+)=
-{1\ov V_{n+1}^+}\,(T_{n+1}(f)\inv\,\dl^+,\,\La'\dl^+).\label{id}\eq

To evaluate the inner product on the right side we for the moment
replace $n+1$ by $n$
so that we can apply earlier formulas. The inner product becomes
\[(\Ti \dl^+,\,\La'\dl^+)=(\La\,\Ti\dlp,\,\dlp).\]
Multiplying the second identity of (\ref{Lcom}) left and right by $\Ti$
gives
\bq\La\,\Ti+\Ti\fp\tn\Tti\dlp=\Ti\,\La+\Ti\dlm\tn\Tti\fm.\label{Ticom}\eq
Applying this to $\dl^+$ (observing that $\La\dlp=0$) and taking the
inner product
with $\dl^+$ we find that
\[(\La\,\Ti\dlp,\,\dlp)=\Vtm_n\,U_n^--U_n^+\,V_n^+.\]
Here we used (\ref{other}).
Replacing $n$ by $n+1$, we see that (\ref{id}) becomes
\[U_n^+=-{\Vtm_{n+1}\ov V_{n+1}^+}\,U_{n+1}^-+U_{n+1}^+.\]
This gives
\bq U_n^+-U_{n+1}^+=\Utm_n\,U_{n+1}^-\label{UU}\eq
by (\ref{UV}).

Next we apply (\ref{ABCD}) to the $n\times n$ matrix
\[\left(\begin{array}{ccccc}0&f_{-1}&\cd&f_{-n}\\
0&f_0&\cd&f_{-n+1}\\
\vdots&\vdots&\cd&\vdots\\
1&f_{n-2}&\cd&f_{-1}\\
0&f_{n-1}&\cd&f_0\end{array}\right).\]
The formula says that the inverse of the upper-left entry
of its inverse equals
\[-(\Ti\,\La\dlm,\,\ftp).\]
This inverse also equals $(-1)^{n-1}\,D_n/K$ where $K$ is the $n-1,\,0$
cofactor of the
matrix (\ref{Tn+1}). But $(T_{n+1}(f)\inv\La\dlm,\,\dlp)$, the
$0,\,n-1$ entry of the inverse of (5), equals $(-1)^{n-1}\,K/D_{n+1}$.
We have shown
\[{D_{n+1}\ov D_n}\,(T_{n+1}(f)\inv\,\La\,\dlm,\,\dlp)
=-(\Ti\,\La\,\dlm,\,\ftp).\]
This may be written
\[(\dlp,\,\La\,T_{n+1}(f)\inv\,\dlm)=-V_{n+1}\,(\dlp,\,\La\,\Ti\,\ftm).\]
To compute the left side we use the fact that
\[\La\,\tl
f^-=\left(\begin{array}{c}f_{-n+1}\\\vdots\\f_{-1}\\0\end{array}\right)=
T_n(f)\,\dlm-f_0\,\dl^-.\]
Hence
\[\Ti\,\La\,\ftm=\dlm-f_0\;\Ti\dlm.\]
Therefore applying (\ref{Ticom}) to $\ftm$ gives
\[\La\,\Ti\,\ftm=\dlm-f_0\;\Ti\dlm\]\[+(\Tti\fm,\,\ftm)\;\Ti\dlm
-(\Tti\,\dlp,\,\ftm)\;\Ti\fp.\]
Taking inner products with $\dlp$ we find
\[(\dlp,\,\La\,\Ti\,\ftm)=-\left(f_0-(\Tti\fm,\,\ftm\right)\,\Vtm_n
-\Utm_n\,\Up_n.\]
But (\ref{fn}) gives
\[f_0-(\Tti\fm,\,\ftm)=f_0-(\Ti\fp,\,\ftp)={1\ov \Vp_{n+1}}.\]
Hence
\[(\dlp,\,\La\,\Ti\,\ftm)=-{\Vtm_n\ov \Vp_{n+1}}-\Utm_n\,\Up_n,\]
\bq(\dlp,\,\La\,T_{n+1}(f)\inv\,\dlm)=
\Vtm_n+\Utm_n\,\Up_n\,\Vp_{n+1}.\label{Ldmdp}\eq\sp

\begin{center} {\bf 2. Nonuniversal recursion relations}\end{center}

Here we restrict to our symbol
\[f_I(z)=e^{t/z}\,(1+z)^k\]
but write $f$ instead of $f_I$ for notational convenience. We shall
obtain relations which follow from the representation
\[f_j={1\ov 2\pi i}\int e^{t/z}\,(1+z)^k\,z^{-j-1}\,dz\]
upon integrating by parts. The fact
\[{1\ov 2\pi i}\int {d\ov
dz}\left\{e^{t/z}\,(1+z)^{k+1}\,z^{-j}\right\}\,dz=0\]
gives
\[-t\,(f_{j+1}+f_j)+(k+1)\,f_{j-1}-j\,(f_j+f_{j-1})=0,\]
\[(j+t)\,f_j+(j-k-1)\,f_{j-1}+t\,f_{j+1}=0.\]
Replacing $j$ by $i-j$ we obtain the $i,j$ entry of a matrix identity:
\[(M+t)\,\T-\T\,M+(M-k-1)\,T_n(zf)-T_n(zf)\,M+t\,T(z\inv f)=0.\]
Here $M$ denotes the diagonal matrix with diagonal entries
$1,\,2,\cd,n$.
We use the identities (\ref{Lcom}) and (\ref{L'com}) to write this as
\[(M+t)\,\T-\T\,M+(M-k-1)\,(\La'\,\T+\dlp\tn\ftp)\]
\[-(\T\,\La'+\ftm\tn\dlm)M+t\,(\T\,\La+\fp\tn\dlp)=0,\]
\[(M+t)\,\T-\T\,M+(1-k)\,\La'\,\T-\T\,\La'+t\,\T\,\La\]
\[-k\,\dlp\tn\ftp-n\,\ftm\tn\dlm+t\,\fp\tn\dlp=0.\]
We multiply this left and right by $\Ti$, obtaining
\[\Ti\,(M+t)-M\,\Ti+\Ti\,(M-k-1)\,\La'-\La'\,M\,\Ti+t\,\La\,\Ti\]
\bq-k\,\Ti\dlp\tn\Tti\ftp-n\,
\Ti\ftm\tn\Tti\dlm+t\,\Ti\fp\tn\Tti\dlp=0.\label{Mid}\eq
This is the basic matrix identity. Applying this matrix identity to
$\dl^{\pm}$ and
taking inner products with $\dl^{\pm}$ gives identities for scalar
quantities. We shall
need three of the four.

First we apply (\ref{Mid}) to $\dlp$ and take the inner product with
$\dlp$. If we
recall our definitions and the fact $\La\dlp=0$ we obtain
\[t\,\Vp_n+(1-k)\,(\Ti\,\La'\,\dlp,\,\dlp)+t\,(\La\,\Ti\dlp,\,\dlp)\]
\bq-k\,\Vp_n\,\Utp_n-n\,(\Ti\ftm,\,\dlp)\,\Vm_n+t\,\Up_n\Vp_n=0.\label{dpdp}\eq
But
\[(\Ti\ftm,\,\dlp)=(\Tti\ftp,\,\dlm)=\Utm_n,\]
so (\ref{UV}) gives
\[\Vm_n=-\Vp_n\,\Um_{n-1}.\]
Thus
\[(\Ti\ftm,\,\dlp)\,\Vm_n=-\Vp_n\,\Utm_n\,\Um_{n-1}.\]
Next, by (\ref{id}) we have
\[(\Ti\,\La'\,\dlp,\,\dlp)=-\Vp_n\,\Utp_{n-1},\ \ \ \
(\La\,\Ti\dlp,\,\dlp)=-\Vp_n\,\Up_{n-1}.\]
Substituting these relations into (\ref{dpdp}) and dividing by $\Vp_n$
gives
\[t+(k-1)\,\Utp_{n-1}-t\,\Up_{n-1}-k\,\Utp_n+n\,\Utm_n\,\Um_{n-1}+t\,\Up_n=0.\]
But (\ref{UU}) gives $\Utm_n\,\Um_{n-1}=\Utp_{n-1}-\Utp_n$, from which
we see
that the above becomes
\[t+(k+n-1)\,\Utp_{n-1}-(k+n)\,\Utp_n+t\,(\Up_n-\Up_{n-1})=0.\]
The derivation holds also for $n=1$ when one defines $\Up_0=\Utp_0=0$,
as
one can easily check. Therefore summing over $n$ gives
\bq n\,t-(k+n)\,\Utp_n+t\,\Up_n=0.\label{id1}\eq

For the next relation we apply (\ref{Mid}) to $\dlm$ and take the inner
product with $\dlp$
to obtain
\[(n+t-1)\,\Vtm_n+t\,(\La\,\Ti\,\dlm,\,\dlp)-k\,\Vp_n\,\Utm_n-n\,\Utm_n\,\Vp_n
+t\,\Up_n\,\Vtm_n=0.\]
If we applying (\ref{Ldmdp}), divide by $\Vp_n$ and use (\ref{UV}) we
obtain
\[-(n+t-1)\,\Utm_{n-1}+t\,\left({\Vtm_{n-1}\ov\Vp_n}+
\Up_{n-1}\,\Utm_{n-1}\right)
-(k+n)\,\Utm_n-t\,\Up_n\,\Utm_{n-1}=0.\]
But
\[{\Vtm_{n-1}\ov\Vp_n}={\Vtm_{n-1}\ov\Vp_{n-1}}\,{\Vtp_{n-1}\ov\Vp_n}
=-\Utm_{n-2}\,(1-\Um_{n-1}\,\Utm_{n-1}),\]
by (\ref{UV}) and (\ref{UV1}). Replacing $n$ by $n+1$ and introducing
the function
\[\ph=1-\Um_n\,\Utm_n\]
we obtain
\[-(t+n)\,\Utm_n+t\,(\Up_n\,\Utm_n-\Utm_{n-1}\,\ph)-(k+n+1)\,\Utm_{n+1}-
t\,\Up_{n+1}\,\Utm_n=0.\]
By (\ref{UU}) this may be written
\bq-(t+n)\,\Utm_n+t\,(\Utm_n)^2\,\Um_{n+1}-t\,
\Utm_{n-1}\,\ph-(k+n+1)\,\Utm_{n+1}=0.
\label{id2}\eq

Finally we we apply (\ref{Mid}) to $\dlm$ and take the inner product
with $\dlm$ to obtain
\[n\,\Vp_n-(n-1)\,(\Ti\dlm,\,\La\,\dlm)-k\,\Vm_n\,
\Utm_n-n\,\Utp_n\,\Vp_n+t\,\Um_n\,\Vtm_n
=0.\]
{}From (\ref{id}) we see that
$(\Ti\dlm,\,\La\,\dlm)=-\Utp_{n-1}\,\Vp_n$.
Substituting
this into the above identity, dividing by $V_n$ and using (\ref{UV})
gives
\bq
t+(n-1)\,\Utp_{n-1}+k\,\Um_{n-1}\,\Utm_n-n\,\Utp_n-t\,\Um_n\,\Utm_{n-1}=0.
\label{id0}\eq
Using (\ref{UU}) we can write
\[k\,\Um_{n-1}\,\Utm_n=n\,\Um_{n-1}\,\Utm_n+(k-n)\,(\Utp_{n-1}-\Utp_n).\]
Substituting this into the preceding and replacing $n$ by $n+1$ give
\bq
t+(k-1)\,\Utp_n+(n+1)\,\Um_n\,\Utm_{n+1}-k\,\Utp_{n+1}-t\,\Um_{n+1}\,\Utm_n=0.
\label{id3}\eq

Another identity can be gotten by using (\ref{UU}) to replace
$n\,(\Utm_{n-1}-\Utm_n)$,
which we see in (\ref{id0}), by $n\,\Um_{n-1}\,\Utm_n$. This gives
\[t-\Utp_{n-1}+(k+n)\,\Um_{n-1}\,\Utm_n-t\,\Um_n\,\Utm_{n-1}=0.\]
Replacing $n$ by $n+1$, multiplying by $\Utm_n$ and adding to
(\ref{id2}) gives
\bq n\,\Utm_n+\Utm_n\,\Utp_n+\ph\,
\left((k+n+1)\,\Utm_{n+1}+t\,\Utm_{n-1}\right)=0.\label{id4}\eq

\begin{center}{\bf 3. Differentiation formulas}\end{center}

We continue to take $f(z)=e^{t/z}\,(1+z)^k$ and write $D_n(t)$ for
$D_n(f)$
Since $df/dt=z\inv f$ we have
\[{d\ov dt}\log\,D_n(t)={\mbox{\rm tr}}\,\Ti\,T_n(z\inv f)=
{\mbox{\rm tr}}\,\left(\La+\Ti\fp\tn\dlp\right)\]
by (\ref{Lcom}), and so
\bq{d\ov dt}\log\,D_n(t)=\Up_n,\label{dlogD}\eq
which is why this quantity arises. Others will arise from further
differentiation.

We use
\[{d\ov dt}\Ti=-\Ti\,T_n(z\inv f)\,\Ti=-\La\,\Ti-\Ti\fp\tn\Tti\dlp,\]\
\[{d\fp\ov dt}=(z\inv f)^+=\La\,\fp+f_{n+1}\,\dlm.\]
Hence
\[{d\Up_n\ov
dt}=-(\La\,\Ti\fp,\,\dlp)-(\Ti\fp,\,\dlp)\;(\Tti\dlp,\,\fp)\]
\[+(\Ti\,\La\,\fp,\,\dlp)+f_{n+1}\,(\Ti\dlm,\,\dlp).\]
By (\ref{Ticom}) the two terms involving $\La$ combine to give
\[(\Ti\fp,\,\dlp)\;(\Tti\dlp,\,\fp)-(\Ti\dlm,\,\dlp)\;(\Tti\fm,\,\fp).\]
Using this we find that the preceding simplifies to
\[-\Vtm_n\;(\Tti\fm,\,\fp)+f_{n+1}\,\Vtm_n=-\Vtm_n\,{1\ov \Vp_{n+1}}\,
{\Vm_{n+2}\ov\Vp_{n+2}},\]
by the second part of (\ref{fn}). This equals
\[-{\Vtm_n\ov\Vp_n}\,{\Vp_n\ov \Vp_{n+1}}\,{\Vtm_{n+2}\ov\Vtp_{n+2}}
=-\Utm_{n-1}\,(1-\Um_n\,\Utm_n)\,\Um_{n+1},\]
by (\ref{UV}) and (\ref{UV1}). We have shown
\bq{d\Up_n\ov dt}=-\ph\,\Utm_{n-1}\,\Um_{n+1}.\label{d1}\eq

In completely analogous fashion (we spare the reader the details) we
compute
\[{d\Um_n\ov dt}=-\Vp_n\;(\Tti\fm,\,\fp)+f_{n+1}\,\Vp_n,\]
and using again the second part of (\ref{fn}), (\ref{UV}) and
(\ref{UV1}) we find that
\bq{d\Um_n\ov dt}=\ph\,\Um_{n+1}.\label{d2}\eq

To find formulas for the derivatives of $\tl U^{\pm}_n$ we use
\[{d\ov dt}\Tti=-\Ti\,T_n(zf)\,\Ti=-\La'\,\Ti+\Ti\ftm\tn\Tti\dlm,\]\
\[{d\ftp\ov dt}=(z\tl f)^+=\La'\,\fp+f_0\,\dlp.\]
At the appropriate points in the computations we use the analogue of
(\ref{Ticom}) with
$\La'$ instead of $\La$, and the first part of (\ref{fn}) rather than
the second.
Again we spare the reader the details. The results are somewhat simpler:
\bq{d\Utp_n\ov dt}=\ph,\label{d3}\eq
\bq{d\Utm_n\ov dt}=-\ph\,\Utm_{n-1}.\label{d4}\eq

Observe that from (\ref{dlogD}), (\ref{id1}), and (\ref{d3}) we have
\[{d\ov dt}t\,{d\ov dt}\log\,D_n(t)=(k+n)\,\ph-n.\]
This gives the representation
\[\log\,D_n(t)=(k+n)\,\int_0^t\log(t/t')\;\ph(t')\,dt'-nt.\]\sp
\pagebreak

\begin{center} {\bf IV. Painlev\'e V and the Laguerre
ensemble}\end{center}
\setcounter{equation}{0}\renewcommand{\theequation}{4.\arabic{equation}}

\begin{center}{\bf 1. Derivation of the differential
equation}\end{center}

We begin by showing how differentiation formulas (\ref{d2}) and
(\ref{d4}) have
analogues in which only indices $n$ and $n-1$ appear on the right side
of (\ref{d2})
and only indices $n$ and $n+1$ appear on the right side of (\ref{d4}).

Solve (\ref{id1})
for $\Um_{n+1}$. The solution involves $\Utm_{n+1}$, which we can solve
for in (\ref{id3}).
Thus $\Um_{n+1}$, and so $d\Um_n/dt$, can be expressed in terms of
quantities with indices
$n$ or $n-1$. To obtain a differentiation formula
for $\Utm_n$ that involves only $n$ and $n+1$ simply solve (\ref{id3})
for $\Utm_{n-1}$.
The results of this are
\bq
{d\Um_n\ov dt}=-{n\ov t}\,\Um_n+{1\ov t}\,\left(\Utp_n-t\ph\right)
{\Um_n\ov \ph -1}+{\ph\ov \ph-1} \Utm_{n-1} {\Um_n}^2,\label{d2Down} \eq
\bq{d\Utm_n\ov dt}={n\ov t}\, \Utm_n +{1\ov t}\, \Utp_n\Utm_n+{1\ov
t}(k+1+n)\,\ph
\,\Utm_{n+1}.\label{d4Up}
\eq

Now compute $\ph'$ using (\ref{d2}) and (\ref{d4Up}). The result can be
put in the form
\bq
{\ph' \ov \ph}-{n\ov t} {\ph-1\ov \ph} - {1\ov t}{\Utp_n(\ph-1)\ov \ph}=
-\Um_{n+1}\Utm_n-{1\ov t}(n+k+1)\Utm_{n+1} \Um_n.
\label{diff1}\eq
Now solve (\ref{id2}) for $\Um_{n+1}\Utm_n$ and insert the result in the
right hand side
of (\ref{diff1}). The result can be written as
\bq k\Utp_{n+1}-(2n+2+k) \Um_n\Utm_{n+1}=t\,{\ph'\ov\ph}
- n {\ph-1\ov \ph}-{\Utp_n(\ph-1)\ov \ph}+t-\Utp_n=\cE,\label{E}\eq
say. Noting that
\bq\Um_n\, \Utm_{n+1}=\Utp_n-\Utp_{n+1}\label{diff4}\eq
we have
\[ \left(\Um_n\,\Utm_{n+1}\right)'=\ph-\Phi_{n+1}\]
by (\ref{d3}). Therefore differentiating the left side of (\ref{E})
gives
\[\cE'=-(2n+2+k)\,\ph+2(k+n+1)\,\Phi_{n+1},\]
\[2(k+n+1)\,\Phi_{n+1}=\cE'+(2n+2+k)\,\ph.\]
Computing $\cE'$ using the right side of (\ref{E}) we obtain the
representation
\bq 2(n+k+1)\,\Phi_{n+1}=2+2(n+k)\,\ph -(n+\Utp_n-\ph){\ph'\ov \ph^2}
- t\,{(\ph')^2\ov \ph^2}+t\,{\ph''\ov\ph}.\label{diff5}\eq
Simply integrating the preceding equation using (\ref{d3}) gives
\bq 2(n+k+1)\,\Utp_{n+1}=t\,{\ph'\ov \ph}-
n\,{\ph-1\ov\ph}+2(k+n)\,\Utp_n+{\Utp_n\ov\ph}+t.
\label{diff6}\eq

We now write out the last relation from which the differential equation
will follow.
\bq
\Phi_{n+1}=1-\Um_{n+1}\Utm_{n+1}
=1-{(\Um_{n+1}\Utm_n)\,(\Um_n\Utm_{n+1})\ov \Um_n\Utm_n}
=1-{(\Um_{n+1}\Utm_n)\,\Um_n\Utm_{n+1})\ov 1-\ph}.\label{diff3}
\eq
Use (\ref{diff4}) and (\ref{diff6}) to express $\Um_n\,\Utm_{n+1}$
in terms of $\Utp_n$, $\ph$ and $\ph'$.  Similarly use (\ref{diff1})
(\ref{diff4}) and (\ref{diff6}) to express $\Utm_n\,\Um_{n+1}$ in
terms of the same quantities.  Finally, on the left hand side of
(\ref{diff3}) use
(\ref{diff5}).  The result is a third-order differential equation
for $\Utp_n$.  (Recall (\ref{d3}).)  This third-order equation is
\[
w^{\prime\prime\prime}={1\ov 2}\left({1\ov w'}+{1\ov w'
-1}\right)(w'')^2
-{1\ov t}\, w'' +{2(k+n)\ov t}\, w' \]
\bq-{2(k+n)\ov t}\, (w')^2 +
{t+n\ov 2 t^2}\,\left(n-t+2\, w\right)-{(n+w)^2\ov 2 t^2 w'}-
{(t-w)^2 \ov 2 t^2 (w' -1)}. \label{de3} \eq

Cosgrove tells us\footnote{Cosgrove,
in his analysis of certain third-order differential equations,
has shown that the third order differential equation of
Chazy Class I (see (A.3) in \cite{cosgrove2}) can be
integrated to a second order and second degree ``master
Painlev\'e equation'' (see (A.21) in \cite{cosgrove2}).
This master Painlev\'e equation, called SD-I in~\cite{cosgrove1},
contains all the Painlev\'e equations I--VI.  Our (\ref{de3})
is a special case of Cosgrove's (A.3).  Carrying out this
reduction~\cite{cosgrove3} in this special case results in
(\ref{firstInt}).  Cosgrove's integration constant equals $-n^2/4$ in
our case.
This  follows from the boundary conditions derived below in (\ref{t0}).}
 that the equation integrates to

\[t^2\,\left(w''\right)^2=-4(k+n)\,t\,\left(w'\right)^3
+\left\{ 4(k+n)\,w+t^2+2(2k+3n)\,t+n^2\right\}\left(w'\right)^2\]
\bq -\left\{2(t+2k+3n)\,w+2n\,t+2n^2\right\}\,w'+(w+n)^2.
\label{firstInt}\eq

The $\s=\s(t)$ form of Painlev\'e V as given by equation (C.45) in
Jimbo-Miwa~\cite{jimbo} (see also \cite{okamoto}) is,
after changing $\s$ to $-\s$ and taking the special parameter values
$\nu_0=\nu_1=0,\ \linebreak \nu_2=k,\ \nu_3=k+n$,
\bq\left(t\,\s''\right)^2=\left\{\s-t\,\s'-2\,(\s')^2+
(2k+n)\,\s'\right\}^2-4\,(\s')^2\,(\s'-k)\,(\s'-k-n).\label{s5}\eq
If
\[w=t-{\s\ov (k+n)},\]
then (\ref{firstInt}) and (\ref{s5}) are equivalent.
Notice that since $w=\Utp_n$, (\ref{id1}) says that $\s=k\,t-t\,\Up_n$
and therefore by
(\ref{dlogD})
\[\s=-t\,{d\ov dt}\log\,\left(e^{-k\,t}\,D_n(t)\right)\]
and therefore
\bq e^{-k\,t}\,D_n(t)
=\exp\left(-\int_0^t{\s(t')\ov t'}\, dt'\right).\label{s5Rep}\eq

This, with Theorem 1, gives Theorem 2 for $G_I$. For $G_D$ it is simply
a matter
of changing $k$ to $-k$ and $t$ to $-t$.

\begin{center} {\bf 2. Laguerre ensemble interpretation of
\mbox{\boldmath $D_n(t)$}}
\end{center}

In order to specify which solution of (\ref{s5}) our $\s$ is, we must
determine the
boundary condition $\s$ satisfies at $t=0$. We have by (\ref{d3}) and
(\ref{UV})
\[{d\s\ov
dt}=(k+n)\,\Um_n\,\Utm_n=(k+n)\,{\Vm_{n+1}\,\Vtm_{n+1}\ov(\Vp_{n+1})^2}.\]
Now $\Vp_{n+1}$ is the upper-left entry of $T_{n+1}(f)\inv$ where,
recall,
$f(z)=e^{t/z}\,(1+z)^k$. As $t\ra0$ this approaches the upper-left entry
of $(I+\La')^{-k}$,
which is clearly equal to 1. Equally clearly, the lower-left entry of
the inverse
has limit ${-k\choose n}$, so that
\[\lim_{t\ra0}\Vm_{n+1}={-k\choose n}.\]
To determine the behavior of $\Vtm_{n+1}$, the the upper-right entry of
$T_{n+1}(f)\inv$,
we write
\bq T_{n+1}(f)=(I+\La')^k+\sum_{p>0,\,q\geq0}{t^p\ov p!}{k\choose
q}\La^{(p-q)},
\label{sum}\eq
where $\La^{(j)}$ denotes $\La^j$ if $j\leq0$ and $(\La')^{-j}$ if
$j<0$.
Factoring out $(I+\La')^k$ and taking the inverse gives
\[ T_{n+1}(f)\inv=\left(I+(I+\La')^{-k}\,
\sum_{p>0,\,q\geq0}{t^p\ov p!}{k\choose
q}\La^{(p-q)}\right)\inv\,(I+\La')^{-k}.\]

If we expand out the inverses we get a sum of products. Each product has
factors of
the form $t^{p_i}\La^{(p_i-q_i)}$ and other factors which are
nonnegative powers of $\La'$. Such a product will have a nonzero
upper-right
entry only if $\sum (p_i+q_i)\geq n$. Therefore, since each $p_i\geq1$,
the lowest
power of $t$ which can occur is $n$. Moreover this power occurs only
when all $q_i=0$
and all the nonegative powers of $\La'$ which occur in the product are
0.
This means that we get the same lowest power of $t$
in the upper-right entry of the inverse if in (\ref{sum}) we replace
$(I+\La')^k$ by $I$
and in the sum we only take the terms with $q=0$. This amounts to
replacing
$T_{n+1}(f)$ by
\[\sum_{p\geq 0}{t^p\ov p!}\,\La^p=e^{t\,\La}.\]
The inverse of this operator is $e^{-t\,\La}$ and the upper-right corner
of this matrix is exactly $(-1)^n\,t^n/n!$. Thus
\[\Vtm_{n+1}={(-t)^n\ov n!}+O(t^{n+1}),\]
as $t\ra0$, and so
\[{d\s\ov dt}=(k+n)\,{-k\choose n}\,{(-t)^n\ov n!}+O(t^{n+1})=
{k\ov n!}\,{n+k\choose n}\,t^{n} + O(t^{n+1}).\]
Since $\s(0)=0$,
\bq \s(t)={k\ov(n+1)!}\,{n+k\choose n}\,t^{n+1}+O(t^{n+2}).\label{t0}\eq

Here is the remarkable fact: the same function $\s$ which satisfies the
equation (\ref{s5})
together with the boundary condition (\ref{t0}) gives a representation
for the Fredholm
determinat which equals the distribution function
for the smallest eigenvalue in the Laguerre ensemble of $k\times k$
matrices associated with the weight function $x^n\,e^{-x}$. Precisely,
we have
\[\mbox{Prob}\,\left(\la_{\mbox{\scriptsize min}}\geq
t\right)=\det\left(I-K_L\right),\]
where $K_L$ is the integral operator on $(0,\,t)$ with kernel
\[K_L(x,y)=\left[k(k+n)\right]^{1/2}\,{\vp_{L,k}(x)\,\vp_{L,k-1}(y)
-\vp_{L,k-1}(x)\,\vp_{L,k}(y)\ov x-y}.\]
Here
\[\vp_{L,k}(x)=\sqrt{{k!\ov (n+k)!}}\,x^{n/2}\,e^{-x/2}\,L^{(n)}_k(x).\]
Moreover
\bq\det\left(I-K_L\right)=\exp\left(-\int_0^t{\s(t')\ov t'}\,
dt'\right).\label{detKLrep}\eq
(See \cite{tw2}, Section VB.)
 It follows from this and (\ref{s5Rep}) that
\[ e^{-k\,t}\,D_n(t)=\det\left(I-K_L\right),\]
which, with Theorem 1, is the assertion of Theorem 3.

\begin{center}{\bf 3. Limiting distribution as
\mbox{$N\ra\iy$}}\end{center}

The foregoing can be restated in more concrete terms as
\bq\sum_{N\ge0}F_I(n;k,\,N)\,{(kt)^N\ov N!}=
e^{kt}\;c_{k,\,n}\;\int_t^{\iy}\cd\int_t^{\iy}\prod x_j^n\;e^{-\sum
x_j}\,\Dl(x)^2\,
dx_1\cd dx_k,\label{Grep}\eq
where $\Dl(x)=\prod_{i<j}(x_i-x_j)$
and $c_{k,\,n}$ is the normalization constant defined by
\[(c_{k,\,n})\inv=\int_0^{\iy}\cd\int_0^{\iy}\prod x_j^n\;e^{-\sum
x_j}\,\Dl(x)^2\;
dx_1\cd dx_k.\]
In fact (\cite{mehta}, formula (17.6.5))
\bq
(c_{k,\,n})\inv=1!\,2!\,\cd\,k!\,\prod_{j=0}^{k-1}(n+j)!.\label{c}\eq

If we make the variable changes $x_j\ra x_j+t$ in the integral, the
right side of
(\ref{Grep}) becomes
\[c_{k,\,n}\;\int_0^{\iy}\cd\int_0^{\iy}\prod (x_j+t)^n\;e^{-\sum
x_j}\,\Dl(x)^2\,
dx_1\cd dx_k.\]
Therefore
\[ F_I(n;k,\,N)={N!\ov
k^N}\,c_{k,\,n}\,\int_0^{\iy}\cd\int_0^{\iy}e^{-\sum x_j}\,
\Dl(x)^2\,dx\;{1\ov 2\pi i}\int t^{-N-1}\prod (x_j+t)^n\,dt,\]
where the inner integral is taken over a contour surrounding $t=0$ and
we write $dx$
for $dx_1\cd dx_k$.

Set
\bq N=kn-r,\ \ \ r=[sk\sqrt {2n}].\label{rdef}\eq
Then
\[{1\ov 2\pi i}\int t^{-N-1}\prod (x_j+t)^n\,dt=
{1\ov 2\pi i}\int
t^r\,\exp\left\{n\,\sum\log\,(1+t\inv\,x_j)\right\}\,{dt\ov t}\]
\[=
{1\ov 2\pi i}\int\exp\left\{r \log t+n\,\left(t\inv\,\sum x_j
-{1\ov 2}t^{-2}\,\sum x_j^2+\cd\right)\right\}\,{dt\ov t},\]
as long as the integration is over a contour where $|t|>\sum x_j$. In
fact we integrate
over the circle
$|t|=n\,\sum x_j/r$. The function $r \log t+n\,t\inv\,\sum x_j$
has a critical point at $t=n\,\sum x_j/r$ and its real part
restricted to the circle has an absolute maximum there. The rest of the
exponent,
$n\left(-{1\ov 2}t^{-2}\,\sum x_j^2+\cd\right)$,
is uniformly bounded on the circle and equals
\[-{r^2\ov 2n}{\sum x_j^2\ov (\sum x_j)^2}+o(1)\]
at the critical point. It follows that as $n\ra\iy$ we have, uniformly
for all $x_j$,
\[{1\ov 2\pi i}\int t^{-N-1}\prod (x_j+t)^n\,dt
\sim \exp\left\{-{r^2\ov2n}\,{\sum x_j^2\ov (\sum x_j)^2}\right\}\,
{1\ov 2\pi i}\int\exp\left\{r \log t+n\,\left({\sum x_j\ov
t}\right)\right\}\,{dt\ov t}\]
\[=\exp\left\{-{r^2\ov2n}\,{\sum x_j^2\ov (\sum
x_j)^2}\right\}\,{(n\,\sum x_j)^r\ov
\Ga(r+1)}.\]
Thus
\bq F_I(n;k,\,N)\sim {N!\ov
k^N}\,c_{k,\,n}\,{n^r\ov\Ga(r+1)}\,\int\cd\int
\exp\left\{-{r^2\ov2n}\,{\sum x_j^2\ov (\sum x_j)^2}\right\}\,(\sum
x_j)^r\,
e^{-\sum x_j}\,\Dl(x)^2\,dx.\label{Frep}\eq

Define
\[\Z:=\{(x_j)\in{\bf R}^k\,:\,\sum x_j=0\}\]
and for general $(x_j)\in{\bf R}^k$ write
\[y=\sum x_j,\ \ \ x_j=-x_j'+y/k,\]
so that $(x_j')\in \Z$. We integrate over $\Z$ with
Lebesgue measure and over $y\in\bf R$. Since each $x_j\ge0$ the $y$
integration
is restricted to
\[y\ge k\,\max x_j'.\]
We find that the double integral in (\ref{Frep}) equals (after changing
back from
$x'$ to $x$)
\[{1\ov\sqrt k}\int_{\Z}\Dl(x)^2\,dx\int_{k\,\max
x_j}^{\iy}\exp\left\{-{r^2\ov2n}\,
{\sum x_j^2+y^2/k\ov y^2}\,\right\}y^r\,e^{-y}\,dy\]
\[={e^{-r^2/2kn}\,\ov\sqrt k}\int_{\Z}\Dl(x)^2\,dx\int_{k\,\max
x_j}^{\iy}
\exp\left\{-{r^2\ov2n}\,{\sum x_j^2\ov y^2}\right\}\,y^r\,e^{-y}\,dy\]
\[={e^{-r^2/2kn}\,\ov\sqrt
k}\,(2n)^{(k^2+r)/2}\,\int_{\Z}\Dl(x)^2\,dx\int_{k\,\max x_j}^{\iy}
\exp\left\{-{r^2\ov2n}\,{\sum x_j^2\ov
y^2}\right\}\,y^r\,e^{-y\,\sqrt{2n}}\,dy,\]
where to obtain the last we made the substitutions
$x_j\ra\sqrt{2n}\,x_j,\
y\ra\sqrt{2n}\,y$.

The factor $y^r\,e^{-y\,\sqrt{2n}}$ in the inner integral achieves its
maximum on
$\bf R^+$ when $y=r/\sqrt{2n}=sk+o(1)$,
at which point the other factor in the integral equals $e^{-\sum
x_j^2}$.
Hence if $\max x_j<s$ (so that $sk$ is interior to the range of the $y$
integration) the
inner intregral is asymptotically equal to
\[e^{-\sum x_j^2}\,\int_0^{\iy}y^r\,e^{-y\,\sqrt{2n}}\,dy=
e^{-\sum x_j^2}\,(2n)^{-(r+1)/2}\,\Ga(r+1),\]
while if $\max x_j>s$ the inner integral is $o$ of this. Moreover the
inner integral is at most
\[\int_{k\,\max x_j}^{\iy}y^r\,e^{-y\,\sqrt{2n}}\,dy
\le e^{-k\,\max x_j}\,\int_0^{\iy}y^r\,e^{-y\,(\sqrt{2n}-1)}\,dy\le C\,
e^{-k\,\max x_j}\,(2n)^{-(r+1)/2}\,\Ga(r+1)\]
for a constant $C$ independent of the $x_j$ and $n$.
Hence application of the dominated convergence theorem shows that if we
define
\[\Z_s=\{x\in\Z:\max x_j\le s\}\]
then the double integral in (\ref{Frep}) is asymptotically
\[{e^{-r^2/2kn}\,\ov\sqrt k}\,(2n)^{(k^2-1)/2}\,\Ga(r+1)\,\int_{\Z_s}
e^{-\sum x_j^2}\,\Dl(x)^2\,dx.\]
If we recall the definition (\ref{rdef}) of $r$ and the value of
$c_{k,\,n}$ given by
(\ref{c}) and apply Stirling's theorem we obtain
\[ F_I(n;k,\,N)\sim \gamma_k\,\int_{\Z_s}e^{-\sum
x_j^2}\,\Dl(x)^2\,dx,\]
where
\bq\ga_k\inv=1!\,2!\,\cd\,k!\;(2\pi)^{(k-1)/2}\;2^{-(k^2-1)/2}.\label{ga}\eq
Equivalently,
\bq\lim_{N\ra\iy}\mbox{Prob}_k\left({\ell_N^I(w)-N/k\ov \sqrt{2N/k}}\le
s\right)
=\ga_k\,\int_{\Z_s}e^{-\sum x_j^2}\,\Dl(x)^2\,dx.\label{Fk0}\eq
The right side is the conditional probability that the largest
eigenvalue
of a matrix from GUE is at most $s$, given that the matrix has trace
zero,
and so Theorem 4 is established.

\begin{center}{\bf 4. Large \mbox{$k$} asymptotics}\end{center}

In this section we denote by $F(s,k)$ the probability in $k\times k$ GUE
that $\la_{\mbox{\scriptsize max}}\le s$ and by $F^0(s,k)$ the same
probability
in $k\times k$ traceless GUE. Thus $F^0(s,\,k)$ is given by the right
side of
(\ref{Fk0}) while
\[F(s,\,k)=c_k\,\int_{-\iy}^s\cd\int_{-\iy}^se^{-\sum
x_j^2}\,\Dl(x)^2\,dx_1\cd dx_k,\]
where (see (\ref{ga}) and formula (17.6.7) of \cite{mehta})
$c_k=\ga_k/\sqrt{\pi}$.
Using now the variable change
\[y=\sum x_j,\ \ \ x_j=x_j'+y/k\]
and then changing $x'$ back to $x$ as before we find that
\[F(s,\,k)={\ga_k\ov\sqrt{\pi k}}\int_{-\iy}^{\iy}dy\,\int_{\Z_{s-y/k}}
e^{-(\sum x_j^2+{y^2\ov k})}\,\Dl(x)^2\,dx\]
\[=\ga_k\,\sqrt{k\ov\pi}\int_{-\iy}^{\iy}e^{-k\,y^2}\,dy\,\int_{\Z_{s-y}}
e^{-\sum
x_j^2}\,\Dl(x)^2\,dx=\sqrt{k\ov\pi}\,\int_{-\iy}^{\iy}e^{-k\,y^2}\,
F^0(s-y,\,k)\,dy.\]
Replacing $s$ by $\sqrt{2k}+s/\sqrt2 k^{1/6}$ and making the variable
change
$y\ra y/\sqrt2 k^{1/6}$ we obtain
\[F(\sqrt{2k}+s/\sqrt2 k^{1/6},\,k)={k^{1/3}\ov\sqrt{2\pi}}\,
\int_{-\iy}^{\iy}e^{-k^{2/3}\,y^2/2}\,
F^0(\sqrt{2k}+(s-y)/\sqrt2 k^{1/6},\,k)\,dy.\]
The factor $F^0(\sqrt{2k}+(s-y)/\sqrt2 k^{1/6},\,k)$ in the integrand is
a bounded
and nonincreasing function of $y$ whereas the integral over $|y|>\dl$ is
$o(1)$ for any $\dl$ and the right side without this factor equals one.
Hence for any
$\dl>0$
\[F^0(\sqrt{2k}+(s-\dl)/\sqrt2 k^{1/6},\,k)+o(1)
\le F(\sqrt{2k}+s/\sqrt2 k^{1/6},\,k)\]
\[\le F^0(\sqrt{2k}+(s+\dl)/\sqrt2 k^{1/6},\,k)+o(1).\]
It follows from this, (\ref{F2}) and the fact that $F_2(s)$ is
continuous that
\bq\lim_{k\ra\iy}F^0(\sqrt{2k}+s/\sqrt2
k^{1/6},\,k)=F_2(s).\label{Fklim}\eq
(A result which includes this, where $N\ra\iy$ and $k\ra\iy$ with $N\gg
k$, can be
found in \cite{johansson3}.)

Another way of expressing (\ref{Fklim}) is as follows. Let $\ell_k$
equal the weak limit of
\linebreak
$(\ell_N^I-N/k)/\sqrt{N}$ as $N\ra\iy$. The
distribution function of $\ell_k$ is
\[\lim_{N\ra\iy}\textrm{Prob}_k\left({\ell_N^I(w)-N/k \ov \sqrt{N}}\le s
\right)
=F^0(s\sqrt{k/2},\,k).\]
Then (\ref{Fklim}) is equivalent to the statement
\[ \lim_{k\ra\iy}\textrm{Prob}\left((\ell_k-2)k^{2/3}\le s\right)=F_2(s)
.\]

\begin{center}{\bf Acknowlegements}\end{center}

The authors have benefited from conversations
and correspondence with D.~Aldous, Y.~Chen, P.~Clarkson,
C.~Cosgrove,  P.~Diaconis, I.~Gessel, C.~Grinstead, A.~Its,
K.~Johansson, J.~L.~Snell and
R.~Stanley.
It is a pleasure to acknowledge this and, as well, the support provided
by the National
Science Foundation through grants DMS--9802122 and DMS--9732687.
Finally, we wish to thank D.~Eisenbud and H.~Rossi for their
support during the MSRI semester Random Matrix Models and Their
Applications.

\end{document}